\documentclass[12pt,a4paper,twoside]{article}

\usepackage[a4paper,margin=1in]{geometry}
\usepackage{amsfonts,amsmath,amsthm,amssymb,mathrsfs}
\usepackage{enumerate}
\usepackage{color,xcolor}
\usepackage{hyperref}
\newcommand{\diam}{\mbox{\rm{diam}}}

\newcommand{\M}{\mathcal}

\allowdisplaybreaks
\numberwithin{equation}{section}

\newtheorem{thm}{Theorem}[section]
\newtheorem{definition}[thm]{Definition}
\newtheorem{prop}[thm]{Proposition}
\newtheorem{lem}[thm]{Lemma}
\newtheorem{exple}[thm]{Example}

\newtheorem{rem}[thm]{Remark}

\DeclareMathOperator{\dist}{dist}

\hypersetup{
	colorlinks=true,
	linkcolor=blue,
	anchorcolor=blue,
	citecolor=blue
}

\begin{document}
	
	\title{The Descriptive Complexity of Continuous Selections with Prescribed Hausdorff dimension of the Graph}
	
	\author{Xiongwei Zhang$^{1}$ \& Megha Pandey$^{2*}$\\
		\\
		\small $^{1}$School of Mathematics and Statistics, Yulin University,\\
		\small Yulin, Shaanxi 719000, P. R. China\\
		\small $^{2}$School of Mathematics, Northwest University,\\
		\small Xi'an, Shaanxi 710127, P. R. China\\
		\\
		\small {\bf Xiongwei Zhang:} zhangxw1019@163.com\\
		\small {\bf Megha Pandey:} meghapandey1071996@gmail.com}
	
	\date{}
	\maketitle
	
	\renewcommand{\thefootnote}{}
	\footnote{*Corresponding author: Megha Pandey.}
	\renewcommand{\thefootnote}{\arabic{footnote}}
	\setcounter{footnote}{0}

	\begin{abstract}
		In this paper, we study the existence of continuous selections of set-valued functions with prescribed Hausdorff dimension of their graphs. We consider Hausdorff-continuous functions defined on the unit interval whose values are nonempty compact subsets of the real line. We prove that, for every $\beta\in[1,2]$, the collection of all such functions admitting a continuous selection whose graph has Hausdorff dimension $\beta$ is $\Sigma^1_1$-complete and hence is not Borel. Further, we show that this result holds even when the set-valued functions have totally disconnected values. Thus, we obtain a complete classification of the descriptive complexity of the existence problem for all possible Hausdorff dimensions of continuous graphs over the unit interval.
	\end{abstract}


	\section{Introduction}\label{introduction}
	
	The study of continuous selections is an important part of the theory of set-valued functions. A continuous selection of a set-valued function $F$ is a continuous function $f$ satisfying $f(x)\in F(x)$ for every point $x$ in the domain. Such selections arise naturally in problems where a set-valued function represents several possible values or solutions at each point, and one wishes to choose a single value that varies continuously with the points of the domain. The theory of continuous selections has applications in different areas of analysis, including optimization, differential inclusions and related problems; see \cite{CastaingValadier1977}. One of the basic questions in this theory is to determine the conditions under which a set-valued function admits a continuous selection. Michael \cite{Michael1956} established a classical result concerning the existence of continuous selections for lower semicontinuous set-valued functions with nonempty closed convex values. Further developments concerning measurable and continuous selections and related properties of multifunctions can be found in \cite{CastaingValadier1977}.
	
	Apart from the existence of continuous selections, it is also natural to investigate their additional properties. In particular, the Hausdorff dimension of the graph of a continuous selection provides information about its geometric structure. In the case of Hausdorff-continuous functions taking nonempty compact values in the real line, the existence of a continuous selection is always guaranteed. Indeed, for such a function $F$, the map $x\mapsto\min F(x)$ is a continuous selection.
	
	
	The Hausdorff dimensions of graphs of continuous functions have been investigated extensively in the literature. It is well known that the Hausdorff dimension of the graph of a continuous real-valued function defined on a nondegenerate compact interval lies between one and two. Fractal interpolation functions, introduced by Barnsley \cite{Barnsley1986}, provide a method for constructing continuous functions whose graphs may have noninteger fractal dimensions. Mauldin and Williams \cite{MauldinWilliams1986} investigated the Hausdorff dimensions of certain graphs, including those associated with Weierstrass-type functions. Later, Shen \cite{Shen2018} established the Hausdorff dimension formula for the graphs of classical Weierstrass functions. In another direction, Fraser and Hyde \cite{FraserHyde2012} proved that the graph of a prevalent continuous function on the unit interval has Hausdorff dimension two. We refer to \cite{Falconer2014,Mattila1995} for the basic notions and properties of Hausdorff dimension and related concepts from fractal geometry.
	
	The construction of continuous functions with prescribed graph dimensions has also been studied in relation to approximation properties and the structure of spaces of continuous functions. Bonilla et al. \cite{Bonilla2021} investigated lineability and spaceability for classes of continuous functions with prescribed Hausdorff and box dimensions of their graphs. More recently, Liu et al. \cite{LiuZhangShi2025} studied algebraic genericity properties of continuous functions with different fractal dimensions. Further results concerning strong spaceability and related properties of continuous functions with prescribed graph dimensions have been obtained in \cite{LiuShiZhang2026,RaposoPellegrino2026}. These studies concern the structure of families of continuous functions having prescribed dimensional properties.
	
	In the context of fractal approximation, Pandey et al. \cite{PandeyAgrawalSom2022} studied dimensional and approximation properties of multivariate $\alpha$-fractal functions. Further, Pandey et al. \cite{PandeySomVerma2024} introduced set-valued $\alpha$-fractal functions and investigated their structural, approximation and dimensional properties. Motivated by dimensional questions in the set-valued setting, we study the existence of continuous selections with prescribed graph dimension from the viewpoint of descriptive set theory.
	
	More precisely, let $F$ be a Hausdorff-continuous compact-valued function and let $\beta$ be a fixed real number. We are interested in determining whether $F$ admits a continuous selection whose graph has Hausdorff dimension $\beta$. The existence of such a selection depends on the values of $F$. For example, consider the constant set-valued function $F(x)=[-1,1]$ for every $x\in[0,1]$. In this case, every continuous function from $[0,1]$ to $[-1,1]$ is a selection of $F$, and hence there exist continuous selections whose graphs have any prescribed Hausdorff dimension between one and two. On the other hand, if $F$ is a constant function whose value is the middle-thirds Cantor set, then every continuous selection of $F$ must be constant. Indeed, the image of a connected set under a continuous function is connected, whereas the Cantor set is totally disconnected. Consequently, every continuous selection in this case has graph dimension one. These examples show that although a continuous selection may exist, it need not have a prescribed Hausdorff graph dimension.
	
	The above observations lead us to investigate the existence of continuous selections with prescribed graph dimension for a general class of set-valued functions. Instead of finding sufficient conditions for an individual set-valued function to admit such a selection, we study the descriptive complexity of the collection of all functions satisfying this condition. Descriptive set theory provides a suitable framework for studying this problem. In particular, analytic sets arise as projections of Borel subsets of products of Polish spaces, whereas $\Sigma^1_1$-complete sets are the analytic sets to which every analytic subset of Baire space can be continuously reduced; see \cite{Kechris1995}.
	
	The descriptive properties of multifunctions and dimension-related sets have also been investigated in different settings. Gregoriades \cite{Gregoriades2011} studied the descriptive complexity of the sets of points of continuity of multifunctions. Mattila and Mauldin \cite{MattilaMauldin1997} established measurability properties of measure and dimension functions on hyperspaces of compact sets. Further, Marcone and Valenti \cite{MarconeValenti2022} studied the descriptive complexity of Salem sets, which are closely related to Hausdorff and Fourier dimensions. In the present paper, we consider a different problem in which the Hausdorff dimension condition is imposed on the graph of a continuous selection. Although the collection of continuous functions having a prescribed Hausdorff graph dimension is Borel, the requirement that at least one continuous selection satisfies this condition leads to a different question. We investigate the descriptive complexity of this existence condition in the space of Hausdorff-continuous compact-valued functions.
	
	Throughout the paper, let $I=[0,1]$ and $J=[-1,1]$. We denote by $\M{K}(J)$ the collection of all nonempty compact subsets of $J$, equipped with the Hausdorff metric $d_H$. Consider the spaces
	\[ X=C(I,\M{K}(J)),\qquad Y=C(I,J), \]
	where $X$ is equipped with the uniform Hausdorff metric
	\[ d_X(F,G)=\sup_{x\in I}d_H(F(x),G(x)), \]
	and $Y$ is equipped with the uniform metric. For $f\in Y$, we denote by
	\[ \Gamma_f=\{(x,f(x)):x\in I\} \]
	the graph of $f$ in the plane. For $F\in X$, let $\operatorname{Sel}(F)$ denote the collection of all continuous selections of $F$. For a fixed $\beta\in\mathbb R$, define
	\begin{equation}\label{A-beta}
		\M{A}_\beta=\{F\in X:\text{there exists }f\in\operatorname{Sel}(F)\text{ such that }\dim_H\Gamma_f=\beta\}.
	\end{equation}
	
	The main objective of this paper is to determine the descriptive complexity of $\M{A}_\beta$ for every possible Hausdorff dimension of the graph of a continuous function. Our main result is stated as follows.

	\begin{thm}\label{main-theorem}
		For every $\beta\in[1,2]$, the set $\M{A}_\beta$ is $\Sigma^1_1$-complete in $X$. In particular, $\M{A}_\beta$ is not Borel. If $\beta\notin[1,2]$, then $\M{A}_\beta$ is empty.
	\end{thm}
	
	The above result shows that although every function in $X$ admits a continuous selection, the existence of a selection with prescribed Hausdorff graph dimension is a non-Borel condition. Further, the result holds even when the set-valued functions have totally disconnected compact values. Thus, our study establishes a connection between continuous selection theory, fractal geometry and descriptive set theory.

	The paper is organized as follows. In Section~\ref{preliminaries}, we recall the basic notions concerning Hausdorff dimension, compact-valued functions and descriptive set theory, together with some preliminary results used throughout the paper. In Section~\ref{positive-dim-section}, we study continuous selections with prescribed graph dimension $1<\beta\leq2$. In Section~\ref{endpoint-section}, we consider the case $\beta=1$ and establish the corresponding completeness result. In Section~\ref{consequences-section}, we obtain the complete classification and discuss some consequences of the constructions. Finally, we conclude the paper in Section~\ref{conclusion-section}.
	
	\section{Preliminaries}\label{preliminaries}
	
	In this section, we recall the definitions and some basic results concerning Hausdorff dimension, continuous selections and descriptive set theory that will be used throughout the paper. We refer to \cite{Falconer2014,Mattila1995} for Hausdorff dimension and to \cite{Kechris1995} for the notions from descriptive set theory.
	
	\begin{definition}\label{hausdorff-definition}
		Let $(Z,d)$ be a metric space and let $E\subset Z$. For $s\geq0$ and $\delta>0$, the $s$-dimensional Hausdorff premeasure of $E$ at scale $\delta$ is defined by
		\[ \M{H}^s_\delta(E)=\inf\left\{\sum_{i=1}^\infty(\diam U_i)^s:E\subset\bigcup_{i=1}^\infty U_i,\ \diam U_i\leq\delta\right\}, \]
		where the infimum is taken over all countable covers of $E$. The $s$-dimensional Hausdorff measure of $E$ is given by
		\[ \M{H}^s(E)=\lim_{\delta\to0}\M{H}^s_\delta(E). \]
		The Hausdorff dimension of $E$ is defined by
		\[ \dim_H E=\inf\{s\geq0:\M{H}^s(E)=0\}. \]
	\end{definition}
	
	Recall that if $s<\dim_HE$, then $\M{H}^s(E)=\infty$, whereas $\M{H}^s(E)=0$ whenever $s>\dim_HE$. In this paper, we study the Hausdorff dimensions of graphs of continuous real-valued functions. For $f\in C(I,J)$, the graph of $f$ is denoted by
	\[ \Gamma_f=\{(x,f(x)):x\in I\}\subset\mathbb R^2. \]
	The Hausdorff dimension of $\Gamma_f$ is computed with respect to the Euclidean metric on $\mathbb R^2$.
	
	We recall that a map $\varphi:E\to Z'$ between metric spaces is called Lipschitz if there exists $L>0$ such that $d'(\varphi(x),\varphi(y))\leq Ld(x,y)$ for every $x,y\in E$. It is called bi-Lipschitz if it is bijective onto its image and both $\varphi$ and its inverse are Lipschitz. We shall use the following properties of Hausdorff dimension.
	
	\begin{lem}\label{dimension-properties}
		Let $E,E_1,E_2,\ldots$ be subsets of a metric space. The following properties hold.
		\begin{enumerate}[(i)]
			\item If $E\subset F$, then $\dim_HE\leq\dim_HF$.
			\item If $\varphi$ is a Lipschitz map on $E$, then $\dim_H\varphi(E)\leq\dim_HE$. In particular, Hausdorff dimension is preserved under bi-Lipschitz maps.
			\item For every countable collection $\{E_j\}_{j\geq1}$,
			\[ \dim_H\left(\bigcup_{j=1}^\infty E_j\right)=\sup_{j\geq1}\dim_HE_j. \]
			\item If $f:I\to J$ is continuous, then $1\leq\dim_H\Gamma_f\leq2$. In particular, the graph of a piecewise linear function on $I$ has Hausdorff dimension one.
		\end{enumerate}
	\end{lem}
	
	The first three properties are standard consequences of Definition~\ref{hausdorff-definition}; see \cite{Falconer2014,Mattila1995}. For (iv), observe that the projection $(x,f(x))\mapsto x$ from $\Gamma_f$ onto $I$ is Lipschitz. Hence, $\dim_H\Gamma_f\geq\dim_HI=1$. The upper bound follows from $\Gamma_f\subset\mathbb R^2$. Moreover, the graph of a piecewise linear function is a finite union of line segments and therefore has Hausdorff dimension one.

	\begin{definition}\label{hausdorff-metric-definition}
		Let $(Z,d)$ be a metric space and let $\M{K}(Z)$ denote the collection of all nonempty compact subsets of $Z$. For $z\in Z$ and a nonempty set $A\subset Z$, put $\dist(z,A)=\inf_{a\in A}d(z,a)$. The Hausdorff distance between $A,B\in\M{K}(Z)$ is defined by
		\[ d_H(A,B)=\max\left\{\sup_{a\in A}\dist(a,B),\sup_{b\in B}\dist(b,A)\right\}. \]
	\end{definition}
	
	\begin{definition}\label{selection-definition}
		A set-valued function $F:I\to\M{K}(J)$ is called Hausdorff-continuous if $d_H(F(x_n),F(x))\to0$ whenever $x_n\to x$ in $I$. A function $f\in C(I,J)$ is called a continuous selection of $F$ if $f(x)\in F(x)$ for every $x\in I$. We denote the collection of all such selections by
		\[ \operatorname{Sel}(F)=\{f\in C(I,J):f(x)\in F(x)\text{ for every }x\in I\}. \]
	\end{definition}
	
	We consider $X=C(I,\M{K}(J))$ with the uniform Hausdorff metric $d_X$ and $Y=C(I,J)$ with the uniform metric, as introduced in Section~\ref{introduction}. For compact $Z$, the hyperspace $(\M{K}(Z),d_H)$ is compact and metrizable. In particular, $X$ and $Y$ are separable and complete with respect to their uniform metrics; see \cite[Sections~4.E--4.F]{Kechris1995}. Since $J$ is a compact interval, every $F\in X$ admits a continuous selection. Indeed, for nonempty compact subsets $A,B\subset J$,
	\[ |\min A-\min B|\leq d_H(A,B). \]
	Thus, $x\mapsto\min F(x)$ belongs to $\operatorname{Sel}(F)$. The problem studied here concerns the additional requirement that a selection has a prescribed Hausdorff dimension of its graph.
	
	\begin{definition}\label{polish-definition}
		A topological space is called Polish if it is separable and admits a compatible complete metric. The Borel subsets of a topological space are the elements of the smallest $\sigma$-algebra containing all open sets. A set is called $G_\delta$ if it is a countable intersection of open sets.
	\end{definition}
	
	Let $\omega=\{0,1,2,\ldots\}$ with the discrete topology. The product space $\omega^\omega$ is called Baire space, and $D=\{0,1\}^{\mathbb N}$, where $\mathbb N=\{1,2,\ldots\}$, is called Cantor space. Both are Polish spaces, and $D$ is compact. We shall use the product topology on these spaces.
	
	\begin{definition}\label{analytic-definition}
		Let $Z$ be a Polish space. A set $A\subset Z$ is called analytic, or $\Sigma^1_1$, if there exist a Polish space $W$ and a Borel subset $B\subset Z\times W$ such that
		\[ A=\{z\in Z:\text{there exists }w\in W\text{ with }(z,w)\in B\}. \]
		Equivalently, $A$ is the projection of a Borel subset of $Z\times W$ onto the first coordinate.
	\end{definition}
	
	Every Borel set in a Polish space is analytic, but an analytic set need not be Borel; see \cite[Chapter~14]{Kechris1995}. For the lower bounds in this paper, we use continuous reductions between subsets of Polish spaces.
	
	\begin{definition}\label{complete-definition}
		Let $A\subset Z$ and $E\subset W$, where $Z$ and $W$ are Polish spaces. We say that $E$ is continuously reducible to $A$ if there exists a continuous map $\varphi:W\to Z$ satisfying
		\[ w\in E\quad\Longleftrightarrow\quad\varphi(w)\in A. \]
		An analytic set $A\subset Z$ is called $\Sigma^1_1$-complete if every analytic subset of Baire space $\omega^\omega$ is continuously reducible to $A$.
	\end{definition}
	
	The following fact will be used to obtain the non-Borel conclusion of the main theorem.

	\begin{lem}\label{non-Borel-lemma}
		Every $\Sigma^1_1$-complete subset of a Polish space is non-Borel.
	\end{lem}
	
	\begin{proof}
		Let $E\subset\omega^\omega$ be an analytic set which is not Borel; see \cite[Theorem $14.2$]{Kechris1995}. Suppose that $A$ is $\Sigma^1_1$-complete and Borel. By Definition~\ref{complete-definition}, there exists a continuous map $\varphi:\omega^\omega\to X$ such that $E=\varphi^{-1}(A)$. Since the continuous inverse image of a Borel set is Borel, $E$ must be Borel, which is a contradiction. Hence, $A$ is not Borel.
	\end{proof}

	We also need the following measurability result concerning Hausdorff dimension.
	
	\begin{lem}[Mattila and Mauldin \cite{MattilaMauldin1997}]\label{dimension-Borel}
		For a compact metric space $Z$, the function $K\mapsto\dim_H K$ from $\M{K}(Z)$ to $[0,+\infty]$ is Borel measurable.
	\end{lem}
	
	In particular, for every fixed $\beta\in\mathbb R$, the collection $\{K\in\M{K}(Z):\dim_HK=\beta\}$ is Borel. We now establish the analytic upper bound for the set defined in \eqref{A-beta}.
	
	\begin{prop}\label{analytic-upper}
		For every fixed $\beta\in\mathbb R$, the set $\M{A}_\beta$ is analytic in $X$.
	\end{prop}
	
	\begin{proof}
		For $F\in X$ and $f\in Y$, put $L(F,f)=\sup_{x\in I}\dist(f(x),F(x))$. For nonempty compact subsets $A,B\subset J$ and $y,z\in J$, we have $|\dist(y,A)-\dist(z,B)|\leq|y-z|+d_H(A,B)$. Consequently,
		\[ |L(F,f)-L(G,g)|\leq\|f-g\|_\infty+d_X(F,G) \]
		for every $F,G\in X$ and $f,g\in Y$. Thus, $(F,f)\mapsto L(F,f)$ is continuous, and $L(F,f)=0$ precisely when $f\in\operatorname{Sel}(F)$. Therefore,
		\[ \M{S}=\{(F,f)\in X\times Y:f\in\operatorname{Sel}(F)\} \]
		is a closed subset of $X\times Y$.
		
		For $f,g\in Y$, pairing points of their graphs with the same first coordinate gives $d_H(\Gamma_f,\Gamma_g)\leq\|f-g\|_\infty$. Hence, the map $f\mapsto\Gamma_f$ from $Y$ to $\M{K}(I\times J)$ is continuous. By Lemma~\ref{dimension-Borel}, the set
		\[ E_\beta=\{f\in Y:\dim_H\Gamma_f=\beta\} \]
		is Borel. Consequently, $\M{S}\cap(X\times E_\beta)$ is Borel, and its projection onto $X$ is exactly $\M{A}_\beta$. This proves the result.
	\end{proof}

	We next recall the complete analytic set that will be used to obtain the lower bounds. Let $D=\{0,1\}^{\mathbb N}$ be Cantor space, where $\mathbb N=\{1,2,\ldots\}$. For two distinct sequences $u,v\in D$, let $m$ denote the length of their longest common initial word and put $d_D(u,v)=2^{-m}$; for $u=v$, take $d_D(u,u)=0$. If $v$ is a finite binary word, we write $[v]$ for the set of all infinite binary sequences beginning with $v$. The space $D$ is compact, and its coordinate projections are continuous maps into the discrete space $\{0,1\}$; see \cite{Kechris1995}.
	
	\noindent
	Define
	\begin{equation}\label{P-infty}
		P_\infty=\{u\in D:u_j=1\text{ for infinitely many }j\},\quad \M{P}=\{K\in\M{K}(D):K\cap P_\infty\neq\varnothing\}.
	\end{equation}

	\begin{lem}\label{analytic-complete-test}
		The set $\M{P}$ is $\Sigma^1_1$-complete in $\M{K}(D)$.
	\end{lem}
	
	\begin{proof}
		Observe that
		\[ P_\infty=\bigcap_{m=1}^\infty\bigcup_{j\geq m}\{u\in D:u_j=1\}. \]
		Thus, $P_\infty$ is a $G_\delta$ subset of $D$. Since the membership relation
		\[ \{(K,u)\in\M{K}(D)\times D:u\in K\} \]
		is closed, the set $\{(K,u):u\in K\cap P_\infty\}$ is Borel. Consequently, its projection $\M{P}$ is analytic.
		
		By \cite[Exercise~27.3 and Section~27.B, p.~209]{Kechris1995}, the collection of compact subsets of Cantor space containing a sequence with infinitely many ones is $\Sigma^1_1$-complete. Since $0^\infty\notin P_\infty$, the same result holds when only nonempty compact subsets are considered, by using the continuous map $K\mapsto K\cup\{0^\infty\}$. Hence, $\M{P}$ is $\Sigma^1_1$-complete.
	\end{proof}

	\section{Prescribed Graph Dimensions Above One}\label{positive-dim-section}
	
	In this section, we prove the $\Sigma^1_1$-completeness of $\M{A}_\beta$ for $1<\beta\leq2$. We first construct continuous functions with prescribed graph dimensions. We will use the properties of Hausdorff dimension recalled in Lemma~\ref{dimension-properties}.
	
	\begin{lem}\label{dimension-profile}
		For every $s\in(1,2)$, there exists a continuous function $w_s:I\to[0,1]$ satisfying
		\[ w_s(0)=w_s(1)=0,\qquad\dim_H\Gamma_{w_s}=s. \]
	\end{lem}
	
	\begin{proof}
		Put $\lambda=2^{s-2}\in(1/2,1)$ and consider the Weierstrass function
		\[ W(t)=\sum_{n=0}^\infty\lambda^n\cos(2\pi 2^nt),\qquad t\in\mathbb R. \]
		The series converges uniformly. By Shen's dimension formula \cite{Shen2018}, the graph of $W$ has Hausdorff dimension $2+\log\lambda/\log2=s$. Since $W$ is one-periodic, the graph over $I$ has the same dimension. Let $M=(1-\lambda)^{-1}$. Then, $|W(t)|\leq M$ and $W(0)=W(1)=M$. Therefore,
		\[ w_s(t)=\frac{M-W(t)}{2M} \]
		has values in $[0,1]$, vanishes at both endpoints, and has a graph of Hausdorff dimension $s$, since its graph is an invertible affine image of the graph of $W|_I$.
	\end{proof}
	
	Fix $\beta\in(1,2]$ and put
	\[ s_j=\beta-(\beta-1)2^{-j},\qquad I_j=[2^{-j},2^{1-j}],\qquad j\geq1. \]
	Then, $1<s_j<\beta$ and $s_j\to\beta$. The intervals $I_j$ have pairwise disjoint interiors and their union is $(0,1]$. Define $h_j:I\to[0,1]$ by
	\begin{equation}\label{h-profile}
		h_j(x)=\begin{cases}w_{s_j}(2^jx-1),&x\in I_j,\\0,&x\notin I_j.\end{cases}
	\end{equation}
	By Lemma~\ref{dimension-profile}, every $h_j$ is continuous.
	
	\begin{lem}\label{positive-family}
		For $u=(u_j)_{j\geq1}\in D$, define
		\begin{equation}\label{positive-g}
			g_u(x)=\sum_{j=1}^\infty4^{-j}u_j(2+h_j(x)),\qquad x\in I.
		\end{equation}
		Then, $(x,u)\mapsto g_u(x)$ is continuous from $I\times D$ to $[0,1]$. If $u$ and $v$ first differ at coordinate $k$, with $u_k=0$ and $v_k=1$, then
		\[ g_v(x)-g_u(x)\geq4^{-k}\qquad\text{for every }x\in I. \]
		Moreover,
		\begin{equation}\label{positive-dimension}
			\dim_H\Gamma_{g_u}=\sup\left(\{1\}\cup\{s_j:u_j=1\}\right).
		\end{equation}
		In particular, $\dim_H\Gamma_{g_u}=\beta$ if and only if $u\in P_\infty$.
	\end{lem}
	
	\begin{proof}
		Since $0\leq h_j(x)\leq1$, the $j$th summand in \eqref{positive-g} lies between zero and $3\cdot4^{-j}$. The series converges uniformly on $I\times D$, and its values belong to $[0,1]$. As each summand is continuous in $(x,u)$, the resulting map is continuous.
		
		Suppose that $u$ and $v$ first differ at coordinate $k$, with $u_k=0$ and $v_k=1$. The contribution of the $k$th coordinate to $g_v-g_u$ is at least $2\cdot4^{-k}$, whereas the total possible negative contribution from the remaining coordinates is at most $3\sum_{j>k}4^{-j}=4^{-k}$. Therefore,
		\[ g_v(x)-g_u(x)\geq2\cdot4^{-k}-3\sum_{j>k}4^{-j}=4^{-k} \]
		for every $x\in I$.
		
		Put $c_u=2\sum_{j=1}^\infty4^{-j}u_j$. Since the functions $h_j$ have disjoint supports apart from their endpoints, where they vanish, we have
		\[ g_u(x)=c_u+4^{-j}u_jh_j(x),\qquad x\in I_j. \]
		If $u_j=0$, this is a constant function on $I_j$ and its graph is a horizontal segment. If $u_j=1$, the graph over $I_j$ is an invertible affine image of $\Gamma_{w_{s_j}}$, and consequently has Hausdorff dimension $s_j$. Since $\Gamma_{g_u}$ consists of these countably many graph pieces together with the point $(0,c_u)$, countable stability of Hausdorff dimension gives \eqref{positive-dimension}. If $u$ has infinitely many ones, the corresponding indices are unbounded, and hence the supremum is $\beta$. If $u$ has only finitely many ones, the supremum is strictly smaller than $\beta$.
	\end{proof}
	
	The following proposition is the main property of the construction. It shows that the continuous selections of the corresponding multifunctions are exactly the functions constructed above.
	
	\begin{prop}\label{family-selections}
		Suppose that $g:I\times D\to J$, $(x,u)\mapsto g_u(x)$, is continuous and, for every two distinct $u,v\in D$, either $g_u(x)<g_v(x)$ for every $x\in I$, or $g_v(x)<g_u(x)$ for every $x\in I$. For $K\in\M{K}(D)$, define
		\begin{equation}\label{F-K}
			F_K(x)=\{g_u(x):u\in K\},\qquad x\in I.
		\end{equation}
		Then, $F_K\in X$, the map $K\mapsto F_K$ is continuous from $\M{K}(D)$ to $X$, and
		\begin{equation}\label{exhaustion}
			\operatorname{Sel}(F_K)=\{g_u:u\in K\}.
		\end{equation}
		Furthermore, $F_K(x)$ is homeomorphic to $K$ for every $x\in I$.
	\end{prop}
	
	\begin{proof}
		Since $g$ is continuous and $K$ is compact, $F_K(x)$ is a nonempty compact subset of $J$ for every $x\in I$. Moreover,
		\[ d_H(F_K(x),F_K(y))\leq\sup_{u\in K}|g_u(x)-g_u(y)|. \]
		By uniform continuity of $g$ on the compact space $I\times D$, the right-hand side tends to zero as $|x-y|\to0$. Thus, $F_K\in X$.
		
		We next show continuity with respect to $K$. Given $\varepsilon>0$, uniform continuity of $g$ gives $\delta>0$ such that $d_D(u,v)<\delta$ implies $\|g_u-g_v\|_\infty<\varepsilon$. If $d_H(K,L)<\delta$, then every $u\in K$ can be paired with some $v\in L$ satisfying $d_D(u,v)<\delta$, and conversely. Hence, $d_X(F_K,F_L)\leq\varepsilon$. Therefore $K\mapsto F_K$ is continuous.
		
		It is clear that each $g_u$ with $u\in K$ is a continuous selection of $F_K$. To prove the reverse inclusion, define
		\[ \Psi_K:I\times K\longrightarrow I\times J,\qquad\Psi_K(x,u)=(x,g_u(x)). \]
		The map $\Psi_K$ is continuous and injective, by the separation assumption. Since $I\times K$ is compact and $I\times J$ is Hausdorff, $\Psi_K$ is a homeomorphism onto its image. Its image is the graph of the set-valued map $F_K$. For any $f\in\operatorname{Sel}(F_K)$, the second coordinate of $\Psi_K^{-1}(x,f(x))$ gives a continuous map $u:I\to K$. Each coordinate map $x\mapsto u_j(x)$ is continuous from the connected interval $I$ into the discrete set $\{0,1\}$, and hence is constant. Therefore $u(x)$ is independent of $x$ and $f=g_u$ for some $u\in K$. This proves \eqref{exhaustion}.
		
		Finally, for fixed $x\in I$, the map $u\mapsto g_u(x)$ is a continuous injection from the compact space $K$ into $J$. It is therefore a homeomorphism of $K$ onto $F_K(x)$.
	\end{proof}
	
	\begin{thm}\label{positive-completeness}
		For every $\beta\in(1,2]$, the set $\M{A}_\beta$ is $\Sigma^1_1$-complete in $X$.
	\end{thm}
	
	\begin{proof}
		Fix $\beta\in(1,2]$ and use the functions $g_u$ given by Lemma~\ref{positive-family}. For $K\in\M{K}(D)$, define $F_K$ as in \eqref{F-K}. By Proposition~\ref{family-selections}, the map $K\mapsto F_K$ is continuous and the selections of $F_K$ are exactly the functions $g_u$ with $u\in K$. Hence, by Lemma~\ref{positive-family},
		\[ F_K\in\M{A}_\beta\ \Longleftrightarrow\ (\exists u\in K)\ \dim_H\Gamma_{g_u}=\beta\ \Longleftrightarrow\ K\cap P_\infty\neq\varnothing. \]
		Thus, $K\mapsto F_K$ is a continuous reduction of the $\Sigma^1_1$-complete set $\M{P}$ to $\M{A}_\beta$. Lemma~\ref{analytic-complete-test} implies that $\M{A}_\beta$ is $\Sigma^1_1$-hard, while Proposition~\ref{analytic-upper} shows that it is analytic. Therefore $\M{A}_\beta$ is $\Sigma^1_1$-complete.
	\end{proof}
	
	\section{Continuous Selections with Graph Dimension One}\label{endpoint-section}
	
	We now consider the endpoint $\beta=1$. The preceding construction uses a sequence of dimensions $s_j>1$ increasing to a prescribed number $\beta>1$, so it does not give the case $\beta=1$. To treat this case, we keep the same set of binary sequences but use a different family of functions. The idea is to approximate a fixed fractal function by piecewise linear functions and to choose between these functions according to the first occurrence of a nonzero coordinate.
	
	Fix a continuous function $w:I\to[0,1]$ such that
	\begin{equation}\label{fractal-w}
		w(0)=w(1)=0,\qquad\dim_H\Gamma_w=\frac32.
	\end{equation}
	Such a function exists by Lemma~\ref{dimension-profile}. By uniform continuity of $w$, for each $n\geq1$ we can choose a piecewise linear function $p_n:I\to[0,1]$ satisfying
	\begin{equation}\label{polygonal-approx}
		p_n(0)=p_n(1)=0,\qquad\|p_n-w\|_\infty\leq4^{-n}.
	\end{equation}
	For example, $p_n$ may be taken to be the polygonal interpolation of $w$ on a sufficiently fine partition of $I$.
	
	\begin{lem}\label{endpoint-phi}
		For $v\in D$, define
		\begin{equation}\label{phi-definition}
			\phi(v)=\begin{cases}p_r,&r=\min\{n\geq1:v_n=1\}\text{ exists},\\w,&v=0^\infty.\end{cases}
		\end{equation}
		Then, $\phi:D\to C(I,[0,1])$ is continuous with respect to the uniform metric. If $v$ and $v'$ first differ at coordinate $k$, then
		\begin{equation}\label{phi-estimate}
			\|\phi(v)-\phi(v')\|_\infty\leq2\cdot4^{-k}.
		\end{equation}
		Moreover,
		\[ \dim_H\Gamma_{\phi(v)}=\begin{cases}1,&v\neq0^\infty,\\3/2,&v=0^\infty.\end{cases} \]
	\end{lem}
	
	\begin{proof}
		Suppose that $v$ and $v'$ first differ at coordinate $k$. If their common initial word contains a one, then the first occurrence of a one is the same for both sequences, so $\phi(v)=\phi(v')$. Otherwise, one of the two functions is $p_k$, while the other is either $p_l$ for some $l>k$ or $w$. By \eqref{polygonal-approx},
		\[ \|p_k-p_l\|_\infty\leq4^{-k}+4^{-l}\leq2\cdot4^{-k},\qquad\|p_k-w\|_\infty\leq4^{-k}. \]
		This proves \eqref{phi-estimate}. In particular, $\phi$ is uniformly continuous on $D$.
		
		If $v\neq0^\infty$, the function $\phi(v)$ is piecewise linear and its graph has Hausdorff dimension one. If $v=0^\infty$, then $\phi(v)=w$ and its graph has Hausdorff dimension $3/2$ by \eqref{fractal-w}.
	\end{proof}
	
	Let $\sigma:D\to D$ denote the left shift defined by $\sigma(u_1,u_2,\ldots)=(u_2,u_3,\ldots)$. For $u\in D$, put $c_u=2\sum_{n=1}^\infty4^{-n}u_n$, so that $0\leq c_u\leq2/3$. On the intervals $I_j=[2^{-j},2^{1-j}]$ introduced in Section~\ref{positive-dim-section}, define
	\begin{equation}\label{endpoint-g}
		g_u(x)=c_u+4^{-(j+1)}\phi(\sigma^ju)(2^jx-1),\qquad x\in I_j,\quad j\geq1,
	\end{equation}
	and put $g_u(0)=c_u$.
	
	\begin{lem}\label{endpoint-family}
		The map $(x,u)\mapsto g_u(x)$ defined by \eqref{endpoint-g} is continuous from $I\times D$ to $[0,1]$. If $u,v\in D$ first differ at coordinate $k$, with $u_k=0$ and $v_k=1$, then
		\begin{equation}\label{endpoint-separation}
			g_v(x)-g_u(x)\geq\frac56\,4^{-k}\qquad\text{for every }x\in I.
		\end{equation}
		Furthermore,
		\begin{equation}\label{endpoint-dimension}
			\dim_H\Gamma_{g_u}=\begin{cases}1,&u\in P_\infty,\\3/2,&u\notin P_\infty.\end{cases}
		\end{equation}
	\end{lem}
	
	\begin{proof}
		Each function $\phi(v)$ is continuous, takes values in $[0,1]$, and vanishes at both endpoints. Therefore, the definitions in \eqref{endpoint-g} agree at the common endpoints of adjacent intervals. Also, for $x\in I_j$,
		\[ |g_u(x)-c_u|\leq4^{-(j+1)}. \]
		The right-hand side tends to zero as $x\to0$, uniformly in $u$, so $g_u$ is continuous at zero. Since $c_u\leq2/3$ and $4^{-(j+1)}\leq1/16$, the values of $g_u$ belong to $[0,1]$.
		
		For $x$ in the interior of any fixed $I_j$, continuity of $(x,u)\mapsto g_u(x)$ follows from Lemma~\ref{endpoint-phi}, continuity of the shift and continuity of the evaluation map $(f,t)\mapsto f(t)$. The same holds at the endpoints of $I_j$ because all the profiles vanish at zero and one. At $x=0$, the uniform bound $|g_u(x)-c_u|\leq4^{-(j+1)}$ and continuity of $u\mapsto c_u$ give joint continuity. Thus, $(x,u)\mapsto g_u(x)$ is continuous on $I\times D$.
		
		Suppose that $u$ and $v$ first differ at coordinate $k$, with $u_k=0$ and $v_k=1$. Since $u_n=v_n$ for $n<k$,
		\begin{equation}\label{baseline-separation}
			c_v-c_u\geq2\cdot4^{-k}-2\sum_{n>k}4^{-n}=\frac43\,4^{-k}.
		\end{equation}
		We next estimate the difference of the perturbation terms. If $j<k$, then $\sigma^ju$ and $\sigma^jv$ first differ at coordinate $k-j$. By \eqref{phi-estimate},
		\[ 4^{-(j+1)}\|\phi(\sigma^ju)-\phi(\sigma^jv)\|_\infty\leq4^{-(j+1)}2\cdot4^{-(k-j)}=\frac12\,4^{-k}. \]
		If $j\geq k$, then $0\leq\phi(\sigma^ju),\phi(\sigma^jv)\leq1$, and hence
		\[ 4^{-(j+1)}\|\phi(\sigma^ju)-\phi(\sigma^jv)\|_\infty\leq4^{-(j+1)}\leq\frac14\,4^{-k}. \]
		Combining these estimates with \eqref{baseline-separation}, we obtain, for every $x\in(0,1]$,
		\[ g_v(x)-g_u(x)\geq\left(\frac43-\frac12\right)4^{-k}=\frac56\,4^{-k}. \]
		At $x=0$, the difference is $c_v-c_u$, so \eqref{endpoint-separation} holds on all of $I$.
		
		It remains to compute the Hausdorff dimension. If $u\in P_\infty$, every shifted sequence $\sigma^ju$ contains a one. By Lemma~\ref{endpoint-phi}, the function $\phi(\sigma^ju)$ is piecewise linear. Therefore, the restriction of $g_u$ to every $I_j$ is piecewise linear, and its graph has Hausdorff dimension one. Since
		\[ \Gamma_{g_u}=\{(0,c_u)\}\cup\bigcup_{j=1}^\infty\{(x,g_u(x)):x\in I_j\}, \]
		countable stability gives $\dim_H\Gamma_{g_u}=1$.
		
		If $u\notin P_\infty$, then $u$ contains only finitely many ones. There exists $N\geq1$ such that $\sigma^ju=0^\infty$ whenever $j\geq N$. For every such $j$, the restriction of $g_u$ to $I_j$ is an invertible affine image of the graph of $w$, and therefore has Hausdorff dimension $3/2$. The restrictions to the remaining intervals are either piecewise linear or affine images of $w$, so their dimensions are at most $3/2$. Countable stability now gives $\dim_H\Gamma_{g_u}=3/2$. This proves \eqref{endpoint-dimension}.
	\end{proof}
	
	\begin{thm}\label{endpoint-completeness}
		The set $\M{A}_1$ is $\Sigma^1_1$-complete in $X$.
	\end{thm}
	
	\begin{proof}
		For $K\in\M{K}(D)$, define
		\[ F_K(x)=\{g_u(x):u\in K\},\qquad x\in I, \]
		where $g_u$ is given by \eqref{endpoint-g}. By Lemma~\ref{endpoint-family}, the family $(x,u)\mapsto g_u(x)$ is continuous and satisfies the separation condition of Proposition~\ref{family-selections}. Consequently, $K\mapsto F_K$ is continuous from $\M{K}(D)$ to $X$, and
		\[ \operatorname{Sel}(F_K)=\{g_u:u\in K\}. \]
		By \eqref{endpoint-dimension},
		\[ F_K\in\M{A}_1\ \Longleftrightarrow\ (\exists u\in K)\ \dim_H\Gamma_{g_u}=1\ \Longleftrightarrow\ K\cap P_\infty\neq\varnothing. \]
		Thus, $K\mapsto F_K$ continuously reduces $\M{P}$ to $\M{A}_1$. By Lemma~\ref{analytic-complete-test}, the set $\M{A}_1$ is $\Sigma^1_1$-hard, while Proposition~\ref{analytic-upper} shows that it is analytic. Hence, $\M{A}_1$ is $\Sigma^1_1$-complete.
	\end{proof}
	
	\section{Consequences and Examples}\label{consequences-section}
	
	We now complete the proof of the main theorem and discuss some consequences of the constructions.
	
	\begin{proof}[Proof of Theorem~\ref{main-theorem}]
		Theorem~\ref{positive-completeness} proves that $\M{A}_\beta$ is $\Sigma^1_1$-complete for every $1<\beta\leq2$, while Theorem~\ref{endpoint-completeness} gives the same result for $\beta=1$. Since every continuous function $f:I\to J$ has a graph whose projection onto $I$ is Lipschitz and surjective, we have $\dim_H\Gamma_f\geq1$. Also, $\Gamma_f\subset I\times J\subset\mathbb R^2$, so $\dim_H\Gamma_f\leq2$. Hence, $\M{A}_\beta=\varnothing$ whenever $\beta\notin[1,2]$. This proves the theorem.
	\end{proof}
	
	\begin{exple}\label{nonattainment-example}
		Let $\beta\in(1,2]$, and consider the family from Lemma~\ref{positive-family}. For $j\geq1$, let $e_j\in D$ be the sequence with one at coordinate $j$ and zero at every other coordinate. Put $K=\{0^\infty\}\cup\{e_j:j\geq1\}$. Since $e_j\to0^\infty$, the set $K$ is compact. By Proposition~\ref{family-selections}, the continuous selections of $F_K$ are exactly $g_{0^\infty},g_{e_1},g_{e_2},\ldots$. Their graph dimensions are $1,s_1,s_2,\ldots$, where $s_j\to\beta$, but none of these dimensions is equal to $\beta$. Thus, a compact family of continuous selections may contain functions with graph dimensions approaching $\beta$ without containing a function of graph dimension $\beta$.
	\end{exple}
	
	\begin{rem}\label{endpoint-observation}
		In the construction of Section~\ref{endpoint-section}, every continuous selection of $F_K$ has graph dimension either $1$ or $3/2$. Nevertheless, the problem of determining whether a selection of dimension one exists is already $\Sigma^1_1$-complete. Thus, the descriptive complexity is not caused by having a large interval of possible graph dimensions within a single family of selections.
	\end{rem}
	
	\begin{rem}\label{fibers-observation}
		By Proposition~\ref{family-selections}, each fiber $F_K(x)$ in both constructions is homeomorphic to the compact subset $K$ of Cantor space. In particular, all these fibers are totally disconnected, and they are contained in $[0,1]$. Therefore, the reductions used to prove Theorem~\ref{main-theorem} already take values among Hausdorff-continuous multifunctions with totally disconnected compact values.
	\end{rem}
	
	
	\section{Conclusion}\label{conclusion-section}
	
	In this paper, we have studied the existence of continuous selections with prescribed Hausdorff dimension of their graphs for Hausdorff-continuous compact-valued functions defined on the unit interval. We have proved that, for every prescribed dimension between one and two, the collection of such set-valued functions admitting a continuous selection with that dimension is $\Sigma^1_1$-complete and hence is not Borel. The proof is based on constructing continuous families of functions indexed by binary sequences and showing that these families describe all continuous selections of the corresponding set-valued functions.
	
	Our results also show that the same descriptive complexity occurs when the values of the constructed set-valued functions are totally disconnected compact sets. In particular, the result for graph dimension one holds even for a family in which every selection has graph dimension either one or $3/2$. Thus, although a continuous selection always exists for the class of set-valued functions considered here, requiring a prescribed Hausdorff dimension of its graph leads to a non-Borel existence problem. This provides a connection between continuous selection theory, the dimensional properties of graphs and descriptive set theory.
	
	\section*{Declarations}
	\noindent{\bf Data availability:} No data were used in the research described in this article.\par
	\medskip
	\noindent{\bf Declaration of competing interest:} The authors declare that they have no competing interests.\par
	\medskip
	\noindent{\bf Authors' contributions:} All authors contributed equally.
	
	\bibliographystyle{abbrv}
	\bibliography{References}
	
\end{document}